\documentclass[11pt,a4paper]{article}
\usepackage[utf8]{inputenc}
\usepackage[T1]{fontenc}
\usepackage{lmodern}
\usepackage[english]{babel}
\usepackage[margin=2.5cm]{geometry}
\usepackage{graphicx}
\usepackage{booktabs}
\usepackage{natbib}
\usepackage{amsmath}
\usepackage{amssymb}
\usepackage{amsthm}
\usepackage{setspace}
\usepackage{caption}
\usepackage{enumitem}
\usepackage{xcolor}
\usepackage{dsfont}
\usepackage{bm}
\usepackage{float}
\usepackage{xr-hyper}
\definecolor{refblue}{RGB}{0,70,160}
\usepackage[colorlinks=true,linkcolor=refblue,citecolor=refblue,
            urlcolor=refblue,filecolor=refblue]{hyperref}
\newcommand{\doi}[1]{doi: \begingroup\urlstyle{rm}%
  \href{https://doi.org/#1}{\nolinkurl{#1}}\endgroup}
\setcitestyle{round}
\let\oldfootnote\footnote
\renewcommand{\footnote}{\fontsize{9}{11}\selectfont\oldfootnote}
\numberwithin{equation}{section}

\title{On Eigenvector Statistics of Large-Dimensional Sample Covariance Matrices with Applications}
\author{
    Weiming Li\thanks{Email: lwm@shufe.edu.cn} \\
    Shanghai University of Finance and Economics
    \and
    Zeng Li\thanks{Email: liz9@sustech.edu.cn} \\
    Department of Statistics and Data Science \\
    Southern University of Science and Technology
    \and
    Qinwen Wang\thanks{Email: wqw@fudan.edu.cn} \\
    School of Data Science, Fudan University
}
\date{}
\newtheorem{thm}{Theorem}[section]
\newtheorem{lemma}[thm]{Lemma}
\newtheorem{remark}[thm]{Remark}
\newtheorem{example}[thm]{Example}

\newtheorem{col}[thm]{Corollary}
\newtheorem{proposition}[thm]{Proposition}

\newcommand{\bS}{{\bf S}}

\newcommand{\bI}{{\bf I}}
\newcommand{\bX}{{\bf X}}

\newcommand{\bx}{{\bf x}}

\newcommand{\bB}{{\bf B}}

\newcommand{\bSig}{{\boldsymbol\Sigma}}

\newcommand{\bu}{{\bf u}}

\newcommand{\bZ}{{\mathbf Z}}

\newcommand{\bv}{{\mathbf v}}

\newcommand{\um}{\underline{m}}

\newcommand{\E}{\mathbb{E}}

\newcommand{\tr}{{\text{\rm tr}}}

\newcommand{\dif}{\mathop{}\!\mathrm{d}}

\allowdisplaybreaks

\newcommand{\rev}[1]{{#1}}
\newenvironment{revision}{\begingroup}{\endgroup}
\begin{document}

\maketitle

\begin{abstract}
We study the eigenvector-weighted spectral statistics
\[
 T_n(\varphi,g)=p^{-1}\tr\{\varphi(\bS_n)g(\bSig)\}
\]
of large-dimensional sample covariance matrices \(\bS_n\) with population
covariance \(\bSig\). These statistics combine functions of sample eigenvalues
with sample--population eigenvector overlaps and arise in statistical losses
that depend on both. As the dimension and sample size grow proportionally,
we establish a joint central limit theorem for
\(\{T_n(\varphi_a,g_a)\}_{a=1}^r\) with
multiple analytic test functions and population weights, centered by
finite-dimensional deterministic equivalents. A key intermediate step is
a joint Gaussian process limit for the weighted resolvent traces
\[
 m_n^g(z)=p^{-1}\tr\{(\bS_n-z\bI_p)^{-1}g(\bSig)\}.
\]
The limiting means and cross-covariances are expressed in terms of the
population spectral distribution. We apply these results in two settings:
one concerns the variance component of minimum-norm prediction risk, and
the other concerns the Frobenius loss of fixed rotationally invariant
covariance estimators. For fixed affine shrinkage
rules, the loss means, variances, and cross-covariances have explicit formulas
in terms of population spectral moments.
\end{abstract}

\textbf{Keywords:} Sample covariance matrix, high-dimensional data, eigenvectors, weighted spectral statistics, prediction risk, covariance shrinkage

\textbf{Mathematics Subject Classification (2020):} 62H10, 62H15
\section{Introduction}
\begingroup
\clubpenalty=10000
\widowpenalty=10000

Sample covariance matrices are fundamental to multivariate analysis and
underpin principal component analysis, covariance estimation, and
prediction. Their eigenvalues describe variation in the data and play a
central role in many spectral estimators. In the high-dimensional regime,
where the dimension is comparable to the sample size, their behavior has
been studied extensively in random matrix theory. The Marchenko--Pastur
(MP) law and its extensions characterize the limiting empirical eigenvalue
distribution \citep{MP67,BSbook}. A central limit theorem (CLT) for
ordinary linear spectral statistics (LSS) was established by
\citet{BS04}, while the fluctuations of extreme eigenvalues have been
studied in \citep{johnstone01,baik2006eigenvalues}.

Nonetheless, many statistical quantities also involve sample eigenvectors. A
rotationally invariant covariance estimator changes the sample eigenvalues
but retains the sample eigenvectors; for a nonspherical population, its
Frobenius loss depends on their alignment with the population eigenspaces.
The variance component of minimum-norm prediction risk has a similar
dependence on the sample eigenspaces \citep{Hastie2022}.
 This motivates a joint fluctuation analysis of sample
eigenvalues and sample--population eigenvector overlaps.
Let \(\bS_n=n^{-1}\bX\bX^*\), where
\(\bX=\bSig^{1/2}\bZ\), \(\bZ\) has independent real standardized
entries, and \(\bSig\) is a deterministic positive definite population
covariance matrix. Here \(^{*}\) denotes the conjugate transpose. We work
in the MP regime, with \(p,n\to\infty\) and
\(p/n\to y\in(0,\infty)\). For a test function \(\varphi\)
and a population weight function \(g\), consider
\[
 T_n(\varphi,g)
 =p^{-1}\tr\{\varphi(\bS_n)g(\bSig)\}
 =p^{-1}\sum_{i=1}^p\varphi(\lambda_i)\bu_i^*g(\bSig)\bu_i,
\]
where \((\lambda_i,\bu_i)\) are the eigenpairs of \(\bS_n\). If
\(\bSig=\sum_{j=1}^pt_j\bv_j\bv_j^*\), then
\[
 \bu_i^*g(\bSig)\bu_i
 =\sum_{j=1}^pg(t_j)|\bu_i^*\bv_j|^2.
\]
Thus \(T_n(\varphi,g)\) combines a function of the sample eigenvalues
with an aggregate of sample--population eigenvector overlaps. The weight
\(g\) selects or reweights population eigenspaces. When \(g\equiv1\),
the statistic reduces to an ordinary LSS. 
Different test-function and weight pairs occur together in a single
statistical loss. For example, for a fixed shrinkage function \(\ell\),
the normalized Frobenius loss satisfies
\[
 p^{-1}\|\ell(\bS_n)-\bSig\|_F^2
 =T_n(\ell^2,1)-2T_n(\ell,t)+p^{-1}\tr(\bSig^2).
\]
Its fluctuations depend on the cross-covariance between the pairs
\((\ell^2,1)\) and \((\ell,t)\). Comparing two shrinkage rules based on
the same sample introduces further mixed terms of this kind.

Eigenvector empirical spectral distributions provide one way to study this
additional information. \citet{BaiMiaoPan2007} established Gaussian limits
for such distributions, and
\citet{PZ08} obtained a CLT for signal-to-interference ratios involving
eigenvector effects. \citet{LedoitPeche2011} derived deterministic
limits of population-weighted resolvents and related them to
eigenvector overlaps and covariance shrinkage. Quantitative convergence
estimates for these distributions were obtained by
\citet{XiaQinBai2013} and \citet{XiYangYin2020}.
At the fluctuation level, the generalized linear spectral statistic
(GLSS) CLTs of \citet{HuYangHan2024} concern
\(\tr\{f(\bS_n)B_n\}\), with a deterministic ancillary matrix
\(B_n\). Their results include \(B_n=g(\bSig)\), and their joint theorem
allows several test functions with a common ancillary matrix. 

While in this paper, we established the joint CLT of \( \{T_n(\varphi_a,g_a)=p^{-1}\tr\{\varphi_a(\bS_n)g_a(\bSig)\}\}_{a=1}^r\) for several  different test-function and weight pairs
\((\varphi_a,g_a)\) computed from the same sample.  Moreover, we considered 
 the weighted resolvent process
\begin{align*}
 m_n^g(z)&=p^{-1}\tr\{(\bS_n-z\bI_p)^{-1}g(\bSig)\}. 
\end{align*}
Our main result establishes joint convergence of these processes for any
fixed finite collection of weights. The limits are jointly Gaussian
analytic processes whose means and cross-covariances
for different weights are expressed as integrals over the population
spectral distribution.   Cauchy integration
then yields a joint CLT for finite collections of analytic test-function
and weight pairs, with fluctuations on the \(p^{-1}\) scale for the
normalized statistics. 

The main technical difficulty is to transfer a joint LSS limit for
perturbed covariance matrices to the derivatives that represent weighted
resolvent traces. We express each trace as a derivative of a log
determinant under a deterministic matrix perturbation and apply the joint
LSS theorem of \citet{LiLiYao2018} to the resulting family. A uniform
second-moment bound on finite-difference errors justifies passing from
the joint log-determinant limit to the derivatives; mixed derivatives
yield the cross-weight covariance kernel. The marginal bias follows from
\citet{HuYangHan2024}. Local resolvent bounds and smooth spectral cutoffs
give tightness on complex contours, including separate contours around
zero and the nonzero spectrum. These steps extend the joint limit from
negative spectral arguments to the analytic tests needed for the
applications.

Our first application is the noise-driven variance component of
minimum-norm prediction risk which involves  \(\sigma^2n^{-1}\tr(\bS_n^+\bSig)\) when \(p>n\). Previous work has described
prediction risk and benign interpolation \citep{Bartlett2020,Hastie2022},
and \citet{li2021asymptotic} established CLTs for minimum-norm prediction
risk under their model assumptions. We obtain a CLT for the variance
component under a general population covariance satisfying our spectral
bounds. 

Our second application concerns the normalized Frobenius loss of
rotationally invariant covariance estimators. For a fixed deterministic
analytic shrinkage function, the loss CLT combines ordinary and
population-weighted spectral statistics and their cross-covariance. For
fixed affine shrinkage rules, we give explicit limiting means, variances,
and cross-covariances of the centered losses in terms of population
spectral moments. These formulas also describe fluctuations of loss
differences. The two-point population
example in Section~\ref{sec:explicit-loss} shows that populations with
the same leading loss for the sample covariance can have different
second-order variability.

The rest of the paper is organized as follows.
Section~\ref{sec:main-results} states the joint limit theorems, and
Section~\ref{sec:applications} develops the two applications.
Section~\ref{sec:proofs} proves the main theorems and application
results, including the two-point example. The supplement contains
supporting estimates and the covariance calculation for the special
case in Theorem~\ref{thm:weight-t}.
\par\endgroup

\section{Main Results}\label{sec:main-results}

\subsection{Preliminaries}\label{sec:preliminaries}

Let \(\bX=(\bx_1,\ldots,\bx_n)\) be a \(p\times n\) data matrix and
\(\bS_n=n^{-1}\bX\bX^*\) its sample covariance matrix. Write
\(\bS_n=\sum_{i=1}^p\lambda_i\bu_i\bu_i^*\). For a test function
\(\varphi\) and a population weight \(g(\bSig)\), define
\[
 T_n(\varphi,g)=p^{-1}\tr\{\varphi(\bS_n)g(\bSig)\}
 =p^{-1}\sum_{i=1}^p\varphi(\lambda_i)\bu_i^*g(\bSig)\bu_i.
\]

The results below concern analytic \(\varphi\). Their proof uses the
resolvent choice \(\varphi_z(x)=(x-z)^{-1}\), which we introduce after
recalling the ordinary spectral distribution.

For any $n\times n$ Hermitian matrix $\bB_n$, its empirical spectral distribution (ESD) is defined by
\[
F^{\bB_n}(x) = \frac{1}{n}\sum_{i=1}^n \mathds{1}_{\{\lambda_i^{\bB_n} \leqslant x\}}.
\]
If $F^{\bB_n}(x)$ converges to a non-random limit $F(x)$ as $n\to\infty$, we call $F(x)$ the limiting spectral distribution (LSD) of $\bB_n$.  Usually we study ESD and LSD through their Stieltjes transform. The Stieltjes transform of any cumulative distribution function $G$ is defined by
\[
m_G(z) = \int \frac{1}{\lambda-z} \dif G(\lambda),\qquad z\in\mathbb{C}^{+}:=\{u+iv,\, u\in\mathbb{R}, v>0\}.
\]

The ESD of \(\bS_n\) is denoted by \(F^{\bS_n}\), and its Stieltjes transform
can be written as
\[
m_{F_n}(z)=\int \frac{1}{x-z} \dif F^{\bS_n}(x)=\frac{1}{p}\tr\left[(\bS_n-z\bI_p)^{-1}\right]
\]

For \(\varphi_z(x)=(x-z)^{-1}\), the eigenvector-weighted statistic becomes
the weighted resolvent
\[
m_n^g(z)=T_n(\varphi_z,g)
=p^{-1}\tr\{(\bS_n-z\bI_p)^{-1}g(\bSig)\}.
\]
When \(g=1\), this reduces to the Stieltjes transform \(m_{F_n}(z)\).

We impose the following assumptions throughout. The population matrix and
its empirical spectral distribution depend on \(p\), although the index on
\(\bSig=\bSig_p\) is suppressed.
\begin{itemize}[leftmargin=1.2cm]
\item[($H_1$)] \(\bX=\bSig^{1/2}\bZ\), where the entries of \(\bZ\) are
independent copies of a real random variable \(z_{11}\), with a distribution
not depending on \(n\), such that
\[
 \E z_{11}=0,\qquad \E z_{11}^2=1,\qquad \E z_{11}^4=3,
 \qquad \E|z_{11}|^{4+\delta}<\infty
\]
for some \(\delta>0\). Gaussianity is not assumed.
\item[($H_2$)] \(\bSig_p\) is deterministic, real symmetric and positive
definite. For fixed constants \(0<h_1\le h_2<\infty\),
\[
 h_1\le\lambda_{\min}(\bSig_p)\le\lambda_{\max}(\bSig_p)\le h_2
 \quad\hbox{for every }p.
\]
Its empirical spectral distribution
\(H_p=p^{-1}\sum_{j=1}^p\delta_{t_j}\) converges weakly to a probability
measure \(H\) on \([h_1,h_2]\).
\item[($H_3$)] \(y_n=p/n\to y\in(0,\infty)\).
\item[($H_4$)] Each weight \(g\) is a fixed, real-valued bounded function on
\([h_1,h_2]\), with finitely many discontinuities, and
\(H(D_g)=0\), where \(D_g\) is its discontinuity set.
For a finite collection of weights this condition is imposed on every
member of the collection.
\end{itemize}
In particular, weak convergence of \(H_p\) then applies to products of
these weights with continuous functions of the population eigenvalues.
The bound in \((H_2)\) concerns the whole population spectrum, including
any eigenvalues whose multiplicities are \(o(p)\).

\begin{lemma}\label{lem:mz}
Under \((H_1)\)--\((H_3)\), the empirical spectral distribution of
\(\bS_n\) converges almost surely to \(F_{y,H}\). Its Stieltjes
transform \(m\) satisfies
\[
 m(z)=\int\frac{dH(t)}{t\{1-y-yzm(z)\}-z}.
\]
The companion Stieltjes transform
\(\um(z)=-(1-y)/z+ym(z)\) is the solution in \(\mathbb C^+\) of
\begin{equation}\label{eq:MP-limit}
 z=-\frac1{\um(z)}+y\int\frac{t}{1+t\um(z)}\,dH(t).
\end{equation}
\end{lemma}

\begin{lemma}\label{lem:mgz}
{Under \((H_1)\)--\((H_4)\),}
\(m_n^g(z)\to m^g(z)\) almost surely for \(z\in\mathbb C^+\), where
\begin{equation*}
 m^g(z)=-\frac1z\int\frac{g(t)}{1+t\um(z)}\,dH(t).
\end{equation*}
\end{lemma}

For the fluctuation results the centering must retain \(y_n\) and
\(H_p\). Let \(\um_{0,n}\) be the companion Stieltjes transform of
\(F_{y_n,H_p}\), defined by the finite-dimensional equation
\begin{equation*}
 z=-\frac1{\um_{0,n}(z)}+
 y_n\int\frac{t}{1+t\um_{0,n}(z)}\,dH_p(t),
 \qquad \um_{0,n}(z)\in\mathbb C^+\quad(z\in\mathbb C^+).
\end{equation*}
Define
\begin{equation}\label{eq:finite-centering}
 m_{0,n}^g(z)=-\frac1z\int\frac{g(t)}{1+t\um_{0,n}(z)}\,dH_p(t),
 \qquad M_n^g(z)=p\{m_n^g(z)-m_{0,n}^g(z)\}.
\end{equation}
Lemmas~\ref{lem:mz} and \ref{lem:mgz} are the usual MP limit and its
population-weighted resolvent form; see \citet{BSbook} and
\citet{LedoitPeche2011}.
These transforms are continued to the lower half-plane by conjugation
and to real intervals outside their spectral supports by analyticity.
No rate of convergence of \(H_p\) or \(y_n\) is required.

{We use the following spectral enclosure.} Choose
\(U>h_2(1+\sqrt y)^2\). If \(y\ne1\), also choose
\(0<L<h_1(1-\sqrt y)^2\), and set
\begin{equation}\label{eq:spectral-enclosure}
 \mathcal K=
 \begin{cases}
 [L,U],&0<y<1,\\
 [0,U],&y=1,\\
 \{0\}\cup[L,U],&y>1.
 \end{cases}
\end{equation}
The finite-dimensional deterministic spectra lie in \(\mathcal K\) for
all sufficiently large \(n\); the sample spectrum lies in
\(\mathcal K\) with probability tending to one. These statements follow
from the extremal eigenvalue bounds for \(n^{-1}\bZ\bZ^*\) and the
bounds in \((H_2)\); see \citet{BSbook}. In the case \(y>1\),
\(L\) bounds the nonzero eigenvalues, not the zero eigenvalues.

Let \(\Omega=\mathbb C\setminus\mathcal K\), and let
\(\mathcal H(\Omega)\) denote the space of analytic functions on
\(\Omega\), with the topology of locally uniform convergence.
The integration contour \(\mathcal C\) is a finite union of positively
oriented, piecewise smooth simple closed curves with mutually disjoint
interiors. Their interiors together contain \(\mathcal K\), and the
curves lie at positive distance from \(\mathcal K\). Thus the winding
number about each point of \(\mathcal K\) is one.

Let \(E_n=\{\operatorname{spec}(\bS_n)\subset\mathcal K\}\).
For all sufficiently large \(n\), define \(M_{n,\mathrm{loc}}^g\) by
\eqref{eq:finite-centering} on \(E_n\) and set it to zero on \(E_n^c\).
It is then a random element of \(\mathcal H(\Omega)\). We set the finitely
many remaining processes to zero as well; this convention has no effect
on any limiting statement.

\subsection{Joint central limit theorems for weighted spectral statistics}

Our goal in this subsection is a joint CLT for the weighted spectral statistics
\(\{T_n(\varphi_a,g_a)\}_{a=1}^r\), with different analytic test functions and
population weights. For \(\varphi\) analytic on a neighborhood of \(\mathcal C\) and its
interior, Cauchy's formula gives, on \(E_n\),
\[
 \varphi(\bS_n)=\frac1{2\pi i}\oint_{\mathcal C}
       \varphi(z)(z\bI_p-\bS_n)^{-1}\,dz.
\]
Recalling \(m_n^g(z)=p^{-1}\tr\{(\bS_n-z\bI_p)^{-1}g(\bSig)\}\),
we obtain
\begin{equation}\label{eq:weighted-lss-cauchy}
 T_n(\varphi,g)=-\frac1{2\pi i}\oint_{\mathcal C}
                   \varphi(z)m_n^g(z)\,dz.
\end{equation}
Define the finite-dimensional deterministic centering by
\begin{equation}\label{eq:weighted-lss-centering}
 \Theta_n(\varphi,g)=-\frac1{2\pi i}\oint_{\mathcal C}
                   \varphi(z)m_{0,n}^g(z)\,dz.
\end{equation}
Since \(M_n^g=p(m_n^g-m_{0,n}^g)\), it follows that
\[
 p\{T_n(\varphi,g)-\Theta_n(\varphi,g)\}
 =-\frac1{2\pi i}\oint_{\mathcal C}\varphi(z)M_n^g(z)\,dz.
\]
Thus a joint functional CLT for the centered weighted resolvents yields
the desired joint CLT by contour integration, including cross-covariances
between different weights. We first define the cross-covariance kernel for
the resolvent limit.
For \(z\in\Omega\), put
\begin{align}
 d_g(z)&=y\{zm^g(z)\}'
       =y\um'(z)\int\frac{tg(t)}{\{1+t\um(z)\}^2}\,dH(t),\notag\\
 q(z,t)&=\frac{t\um(z)}{1+t\um(z)},\qquad
 w_g(z,t)=\frac{t\{\um(z)g(t)-d_g(z)\}}
                   {z\{1+t\um(z)\}^2}.
       \label{eq:q-w}
\end{align}
Any removable singularities in these expressions are filled in by
analytic continuation. With \(q_a(t)=q(z_a,t)\), define
\begin{equation*}
 D(z_1,z_2)=D_{12}=1-y\int q_1(t)q_2(t)\,dH(t).
\end{equation*}
In particular,
\[
 D_{12}=\um(z_1)\um(z_2)
        \frac{z_1-z_2}{\um(z_1)-\um(z_2)},\qquad
 D_{11}=\frac{\um(z_1)^2}{\um'(z_1)},
\]
where diagonal values are obtained by continuity. The stable branch of
\eqref{eq:MP-limit} makes these denominators nonzero on the contours under
consideration.

For two admissible weights \(g,h\), the cross-covariance function is
\begin{align}
 \Gamma_{g,h}(z_1,z_2)
 &=\frac{2y}{D_{12}}\int w_g(z_1,t)w_h(z_2,t)\,dH(t)
 \label{eq:cross-kernel}\\[-2pt]
 &\quad+\frac{2y^2}{D_{12}^2}
 \left\{\int w_g(z_1,t)q_2(t)\,dH(t)\right\}
 \left\{\int q_1(t)w_h(z_2,t)\,dH(t)\right\}.\nonumber
\end{align}
{Write \(\Gamma_g=\Gamma_{g,g}\). This function satisfies
\(\Gamma_{g,h}(z_1,z_2)=\Gamma_{h,g}(z_2,z_1)\), where \(g\) is paired with
\(z_1\) and \(h\) with \(z_2\).}

\begin{thm}[Weighted resolvent CLT]\label{thm:main}

Assume \((H_1)\)--\((H_4)\). For every fixed finite collection of
weights \(g_1,\ldots,g_r\),
\[
 (M_{n,\mathrm{loc}}^{g_1},\ldots,M_{n,\mathrm{loc}}^{g_r})
 \ \Rightarrow\ (M^{g_1},\ldots,M^{g_r})
 \quad\text{in }\mathcal H(\Omega)^r,
\]
where the limit is a jointly complex Gaussian analytic process.
For each weight \(g\), its mean is
\begin{align}
 \mu_g(z)=\E M^g(z)
 &=\int\frac{y\um'(z)\{d_g(z)-\um(z)g(t)\}t^2}
             {z\um(z)\{1+t\um(z)\}^3}\,dH(t)
             \label{eq:lemMean}\\
 &=-\frac{y\int q(z,t)w_g(z,t)\,dH(t)}{D(z,z)}.\nonumber
\end{align}
For \(z_1,z_2\in\Omega\), its bilinear covariance is
\begin{equation*}
 \operatorname{Cov}\{M^{g_a}(z_1),M^{g_b}(z_2)\}
 =\Gamma_{g_a,g_b}(z_1,z_2).
\end{equation*}
Here covariance is bilinear, without conjugating the second variable.
The Hermitian covariance is
\(\E[(M^{g_a}(z_1)-\mu_{g_a}(z_1))
\overline{(M^{g_b}(z_2)-\mu_{g_b}(z_2))}]
=\Gamma_{g_a,g_b}(z_1,\overline{z_2})\).
In particular, restriction to \(\mathcal C\) gives joint weak convergence
in \(C(\mathcal C;\mathbb C)^r\), equipped with the uniform norm.

\end{thm}

\begin{thm}\label{thm:weight-t}

Under the assumptions of Theorem~\ref{thm:main}, take \(g(t)=t\).
Then the same conclusions hold, with the centering \(m_{0,n}^{t}\).
Writing \(\um_a=\um(z_a)\), the covariance has the equivalent form
\begin{align*}
 \Gamma_t(z_1,z_2)
 &=-\frac{2\{\um_1(1-\um_2'/\um_2^2)
                -\um_2(1-\um_1'/\um_1^2)\}}
 {(z_1-z_2)z_1z_2\um_1^2\um_2^2}\\
 &\quad-\frac{2}{z_1z_2(z_1-z_2)^2\um_1^2\um_2^2}
 +\frac{2\um_1'\um_2'}
 {z_1z_2\um_1^2\um_2^2(\um_1-\um_2)^2}.
\end{align*}
The mean is \eqref{eq:lemMean} with \(g(t)=t\).
The values at \(z_1=z_2\) are defined by continuity; the integral form
\eqref{eq:cross-kernel} avoids the cancelling divided differences.

\end{thm}

For the contour representation \eqref{eq:weighted-lss-cauchy},
\rev{if \(\varphi\) is already defined on \([0,\infty)\), we use that
definition in \(\varphi(\bS_n)\), as for the polynomial tests below.
Otherwise, choose a finite-valued Borel extension of its restriction to
\(\mathcal K\), real-valued when \(\varphi\) is real on \(\mathcal K\).
All such extensions agree on \(E_n\), so the weak limit is independent
of this choice.}

The centering \eqref{eq:weighted-lss-centering} retains \(y_n\) and \(H_p\).
Its first-order limit is
\begin{equation*}
 \Theta(\varphi,g)=-\frac1{2\pi i}\oint_{\mathcal C}
                   \varphi(z)m^g(z)\,dz
\end{equation*}
{is the first-order limit of \(\Theta_n(\varphi,g)\).
The fluctuation theorem below uses the finite-dimensional center
\eqref{eq:weighted-lss-centering}.}

\begin{thm}[CLT for weighted spectral statistics]\label{thm:ewss-clt}

Under the assumptions of Theorem~\ref{thm:main}, let
\(\varphi_a\) be analytic on a neighborhood of \(\mathcal C\) and its
interior, and real-valued on \(\mathcal K\), for \(a=1,\ldots,r\).
Then
\[
 \left(p\{T_n(\varphi_a,g_a)-\Theta_n(\varphi_a,g_a)\}\right)_{a=1}^r
 \ \Rightarrow\ N_r(\boldsymbol\mu,\boldsymbol V),
\]
where
\begin{align}
 \mu_a&=-\frac1{2\pi i}\oint_{\mathcal C}
          \varphi_a(z)\mu_{g_a}(z)\,dz,
          \label{eq:lss-mean}\\
 V_{ab}&=\frac1{(2\pi i)^2}\oint_{\mathcal C}\oint_{\mathcal C}
       \varphi_a(z_1)\varphi_b(z_2)
       \Gamma_{g_a,g_b}(z_1,z_2)\,dz_1\,dz_2.
       \label{eq:lss-cross}
\end{align}
The normal distribution is allowed to be degenerate. In the scalar case
write \(\mu(\varphi,g)=\mu_1\) and \(v(\varphi,g)=V_{11}\).

\end{thm}

\begin{remark}
For \(y>1\), the contours around \(\{0\}\) and \([L,U]\) can be chosen
separately. The test function may be defined analytically on the two
separate enclosed neighborhoods; if it vanishes near zero, the zero
contour contributes nothing. More general separated spectral enclosures
can be used provided that both the finite-dimensional deterministic
support and the sample spectrum are contained in them, the latter with
probability tending to one, and the contours stay a fixed positive
distance from these enclosures. {For such an enclosure, spectral confinement is an additional
condition beyond the existence of a gap in the limiting distribution.}
\end{remark}

\begin{revision}
For \(\varphi\equiv1\),
\[
 T_n(1,g)=p^{-1}\tr g(\bSig)=\Theta_n(1,g)
\]
exactly, so the centered statistic vanishes. For \(\varphi(x)=x\),
\(\Theta_n(x,g)=p^{-1}\tr\{\bSig g(\bSig)\}\), and
\begin{align*}
 \E\tr\{\bS_ng(\bSig)\}&=\tr\{\bSig g(\bSig)\},\\
 \operatorname{Cov}\bigl(\tr\{\bS_ng(\bSig)\},\tr\{\bS_nh(\bSig)\}\bigr)
 &=\frac2n\tr\{\bSig^2g(\bSig)h(\bSig)\}.
\end{align*}

\begin{col}[Linear test functions]\label{col:linear-tests}
Under \((H_1)\)--\((H_4)\), let \(g_1,\ldots,g_r\) be fixed admissible
weights. Then
\[
 \left(\tr\{\bS_ng_a(\bSig)\}-\tr\{\bSig g_a(\bSig)\}\right)_{a=1}^r
 \ \Rightarrow\ N_r(\boldsymbol 0,\boldsymbol V),\qquad
 V_{ab}=2y\int t^2g_a(t)g_b(t)\,dH(t).
\]
The limit may be degenerate.
\end{col}
\end{revision}

\section{Applications}\label{sec:applications}

We apply the weighted spectral statistics CLT to two statistical
problems.

\subsection{Variance component of minimum-norm prediction risk}\label{sec:prediction-variance}

Consider the high-dimensional regression model in the over-parametrized
case {\citep{Hastie2022}},
\begin{equation*}
 y_i=\boldsymbol\beta^*\bx_i+\varepsilon_i,\qquad i=1,\ldots,n.
\end{equation*}

Write \(\mathsf X=(\bx_1,\ldots,\bx_n)^*\),
\(\boldsymbol y=(y_1,\ldots,y_n)^*\), and
\(\boldsymbol\varepsilon=(\varepsilon_1,\ldots,\varepsilon_n)^*\).
Assume that the noise vector is independent of the design, with
\(\E\boldsymbol\varepsilon=0\) and
\(\operatorname{Cov}(\boldsymbol\varepsilon)=\sigma^2\bI_n\), where
\(\sigma^2>0\) is fixed. Define the
minimum-norm estimator
\begin{equation*}
 \widehat{\boldsymbol\beta}=(\mathsf X^*\mathsf X)^+
                              \mathsf X^*\boldsymbol y.
\end{equation*}
For an independent test covariate \(\bx_0\) with mean zero and covariance
\(\bSig\), the conditional prediction risk is
\[
 R_{\mathsf X,\boldsymbol\beta}
 =\E\bigl[\{\bx_0^*(\widehat{\boldsymbol\beta}
                     -\boldsymbol\beta)\}^2\mid\mathsf X\bigr].
\]
{The conditional expectation averages over the training noise
and the test covariate. The risk measures prediction of the conditional
mean at a new covariate.}
It decomposes as \(R_{\mathsf X,\boldsymbol\beta}=B_{\mathsf X,\boldsymbol\beta}
+V_{\mathsf X,\boldsymbol\beta}\), where
\begin{align*}
 B_{\mathsf X,\boldsymbol\beta}
 &=\boldsymbol\beta^*\boldsymbol\Pi\bSig\boldsymbol\Pi\boldsymbol\beta,
 &V_{\mathsf X,\boldsymbol\beta}&=\frac{\sigma^2}{n}\tr(\bS_n^+\bSig),\\
 \boldsymbol\Pi&=\bI_p-\bS_n^+\bS_n,
 &\bS_n&=\mathsf X^*\mathsf X/n.
\end{align*}

\rev{The bias is a quadratic form involving the random projection
\(\boldsymbol\Pi\). The variance is a weighted spectral statistic whose
test function is supported on the nonzero spectrum.}

Assume \(y_n\to y>1\), and use the enclosure
\(\mathcal K=\{0\}\cup[L,U]\) in \eqref{eq:spectral-enclosure}.
Let \(\mathcal C_+\) surround \([L,U]\), with zero outside it, and take
\(\varphi_V(z)=z^{-1}\) near \([L,U]\) and \(\varphi_V(z)=0\) near
zero. This is an analytic function on the union of the two separate
neighborhoods. Put
\[
 \Theta_{V,n}=\Theta_n(\varphi_V,t)
       =-\frac1{2\pi i}\oint_{\mathcal C_+}z^{-1}m_{0,n}^{t}(z)\,dz.
\]

\begin{thm}\label{thm:variance-component}
Under the assumptions above,

\[
 p\{V_{\mathsf X,\boldsymbol\beta}-\sigma^2y_n\Theta_{V,n}\}
 \Rightarrow N(\mu_V,\tau_V^2),
\]
where
\begin{align*}
 \mu_V&=-\frac{\sigma^2y}{2\pi i}\oint_{\mathcal C_+}
                         z^{-1}\mu_t(z)\,dz,\\
 \tau_V^2&=\frac{\sigma^4y^2}{(2\pi i)^2}
 \oint_{\mathcal C_+}\oint_{\mathcal C_+}
 z_1^{-1}z_2^{-1}\Gamma_t(z_1,z_2)\,dz_1\,dz_2.
\end{align*}

\end{thm}

\begin{remark}
{The center \(\sigma^2y_n\Theta_{V,n}\) retains both
\(y_n\) and \(H_p\) through the finite-dimensional MP equation.}
A CLT for total prediction risk additionally requires control of the bias
component, whose fluctuations may depend on the orientation of
\(\boldsymbol\beta\) relative to the eigenspaces of \(\bSig\).
\end{remark}

\begin{remark}
When \(\bSig=\bI_p\) and \(y>1\), the leading-order variance limit is
\[
 V_{\mathsf X,\boldsymbol\beta}\longrightarrow\frac{\sigma^2}{y-1}.
\]
{In fact, the finite-dimensional deterministic centering satisfies
\(\Theta_{V,n}=\{y_n(y_n-1)\}^{-1}\) for all sufficiently large \(n\),
so \(\sigma^2y_n\Theta_{V,n}=\sigma^2/(y_n-1)\).}
\end{remark}

\subsection{Frobenius loss of rotationally invariant covariance estimators}\label{sec:explicit-loss}

Many covariance estimators keep the sample eigenvectors and shrink only
the sample eigenvalues \citep{LedoitPeche2011}. Let
{\(\ell\) be a fixed deterministic shrinkage function, real-valued on
\(\mathcal K\) and analytic on a neighborhood of the chosen contours
and their interiors,} and define
\[
 \widehat{\bSig}_{\ell}=\ell(\bS_n)
             =\sum_{i=1}^p\ell(\lambda_i)\bu_i\bu_i^*.
\]
We study the normalized Frobenius loss
\[
 L_n(\ell)=p^{-1}\|\widehat{\bSig}_{\ell}-\bSig\|_F^2.
\]
Expanding the square gives
\begin{equation}\label{eq:frobenius-loss-expansion}
 L_n(\ell)=p^{-1}\tr\{\ell^2(\bS_n)\}
 -2p^{-1}\tr\{\ell(\bS_n)\bSig\}+p^{-1}\tr(\bSig^2).
\end{equation}
The first term is an ordinary linear spectral statistic. The second term
also involves the population weight and is not determined by the sample
eigenvalues alone.

Set
\[
 \Theta_{0,n}(\ell)=\Theta_n(\ell^2,1),\qquad
 \Theta_{1,n}(\ell)=\Theta_n(\ell,t),
\]
and define the finite-dimensional deterministic centering
\begin{equation}\label{eq:loss-centering}
 L_{0,n}(\ell)=\Theta_{0,n}(\ell)-2\Theta_{1,n}(\ell)
                 +p^{-1}\tr(\bSig^2).
\end{equation}

\begin{thm}\label{thm:frobenius-loss}

Under the assumptions of Theorem~\ref{thm:main} and the condition on
\(\ell\) above,
\[
 p\{L_n(\ell)-L_{0,n}(\ell)\}\Rightarrow N(\mu_L,\tau_L^2).
\]
Let \((\varphi_0,g_0)=(\ell^2,1)\) and
\((\varphi_1,g_1)=(\ell,t)\), and put
\[
 \mu_a=\mu(\varphi_a,g_a),\qquad
 v_{ab}=\frac1{(2\pi i)^2}\oint_{\mathcal C}\oint_{\mathcal C}
  \varphi_a(z_1)\varphi_b(z_2)\Gamma_{g_a,g_b}(z_1,z_2)\,dz_1\,dz_2.
\]
Then
\[
 \mu_L=\mu_0-2\mu_1,\qquad
 \tau_L^2=v_{00}+4v_{11}-4v_{01}.
\]
{The mixed term \(v_{01}\) is evaluated using
\(\Gamma_{1,t}\) in \eqref{eq:cross-kernel}.}

\end{thm}

\begin{remark}
The shrinkage function \(\ell\) is deterministic. If it is estimated
from the same data, its estimation error may contribute at the
\(p^{-1}\) scale and must be included in the loss expansion.
\end{remark}

\begin{remark}
The same argument applies to rotationally invariant precision estimators.
If \(\widehat{\boldsymbol\Omega}_{\ell}=\ell(\bS_n)\) estimates
\(\bSig^{-1}\), then
\[
 p^{-1}\|\widehat{\boldsymbol\Omega}_{\ell}-\bSig^{-1}\|_F^2
 =p^{-1}\tr\{\ell^2(\bS_n)\}
 -2p^{-1}\tr\{\ell(\bS_n)\bSig^{-1}\}+p^{-1}\tr(\bSig^{-2}).
\]
{This uses the weights \(1\) and \(t^{-1}\), which are admissible
under \((H_2)\), and the corresponding finite-dimensional center.
The analyticity condition on \(\ell\) is unchanged. When an inverse
function is used for \(y>1\), it is defined on the nonzero component,
with a separate definition near zero as in Section~\ref{sec:prediction-variance}.}
\end{remark}

\begin{revision}
We now specialize Theorem~\ref{thm:frobenius-loss} to the sample
covariance and fixed affine shrinkage estimators. Write
\[
 H_{k,n}=p^{-1}\tr(\bSig^k),\qquad H_k=\int t^k\,dH(t).
\]

\begin{col}\label{cor:sample-loss}
Under \((H_1)\)--\((H_3)\),
\begin{equation}\label{eq:sample-loss-clt}
 p\left\{p^{-1}\|\bS_n-\bSig\|_F^2-y_nH_{1,n}^2\right\}
 \ \Rightarrow\ N\left(yH_2,\,4y^2H_2^2+8y^3H_1^2H_2\right).
\end{equation}
\end{col}
The finite-dimensional center follows from the MP second moment
\(\Theta_n(x^2,1)=H_{2,n}+y_nH_{1,n}^2\) and the weighted first
moment \(\Theta_n(x,t)=H_{2,n}\). In fact,
\[
 \E\{p^{-1}\|\bS_n-\bSig\|_F^2\}
       =y_nH_{1,n}^2+\frac{H_{2,n}}n.
\]
The difference between the finite-sample expectation and the MP center
is therefore \(H_{2,n}/n\). Multiplication by \(p\) gives the limiting
mean \(yH_2\) in \eqref{eq:sample-loss-clt}. At \(\bSig=I\), the
limiting variance is \(4y^2(1+2y)\).

A fixed affine shrinkage rule has the form
\[
 \ell_{a,b}(x)=ax+b,\qquad
 \widehat{\bSig}_{a,b}=a\bS_n+bI,
\]
where \(a,b\) are deterministic constants. In particular,
\(a=1-\rho\) and \(b=\rho\tau\), with \(0\le\rho\le1\) and
\(\tau>0\), shrink the sample covariance toward \(\tau I\).
Define
\begin{equation}\label{eq:affine-center}
 L_{0,n}(a,b)=(a-1)^2H_{2,n}+2b(a-1)H_{1,n}+b^2
                              +a^2y_nH_{1,n}^2.
\end{equation}
The first three terms are \(p^{-1}\|(a-1)\bSig+bI\|_F^2\);
the last is the leading sampling contribution.

\begin{proposition}\label{prop:affine-loss}
Under \((H_1)\)--\((H_3)\), for fixed real \(a,b\),
\[
 p\left\{p^{-1}\|a\bS_n+bI-\bSig\|_F^2-L_{0,n}(a,b)\right\}
 \ \Rightarrow\ N\{a^2yH_2,\,\tau^2(a,b)\},
\]
where
\begin{equation}\label{eq:affine-variance}
 \tau^2(a,b)=4a^4y^2H_2^2+
 8a^2y\int t^2\{(a-1)t+b+ayH_1\}^2\,dH(t).
\end{equation}
For any fixed finite collection of rules, their centered losses converge
jointly. For two pairs \((a,b)\) and \((c,d)\), the limiting
cross-covariance is
\begin{align}
 \mathcal V\bigl((a,b),(c,d)\bigr)
 &=4a^2c^2y^2H_2^2\label{eq:affine-cross}\\
 &\quad+8acy\int t^2\{(a-1)t+b+ayH_1\}
                        \{(c-1)t+d+cyH_1\}\,dH(t).\nonumber
\end{align}
\end{proposition}
The variance is nonnegative and depends only on \(H_1,\ldots,H_4\).
At \(a=1,b=0\), it reduces to
Corollary~\ref{cor:sample-loss}; at \(a=0\), the loss is deterministic.
The joint formula also gives the fluctuation of a difference of losses,
with variance
\(\tau^2(a,b)+\tau^2(c,d)-2\mathcal V((a,b),(c,d))\).
As in Theorem~\ref{thm:frobenius-loss}, the parameters are fixed in advance.

\begin{example}\label{ex:two-point-shrinkage}
Let the population eigenvalues be \(1-\kappa\) and \(1+\kappa\)
in equal proportions, with fixed \(0<\kappa<1\) and even \(p\).
Then \(H_1=H_{1,n}=1\), \(H_2=H_{2,n}=1+\kappa^2\), and
\(\kappa^2\) measures spectral dispersion at a fixed average eigenvalue.
For
\[
 \widehat{\bSig}_a=a\bS_n+(1-a)I,\qquad 0\le a\le1,
\]
the exact finite-dimensional center and the limiting bias are
\[
 L_{0,n}(a,1-a)=(1-a)^2\kappa^2+a^2y_n,
 \qquad \mu_a=a^2y(1+\kappa^2).
\]
Proposition~\ref{prop:affine-loss} gives
\begin{align}
 \tau_a^2
 &=4a^4y^2(1+\kappa^2)^2\label{eq:two-point-variance}\\
 &\quad+4a^2y\bigl[(1-\kappa)^2\{ay+(1-a)\kappa\}^2
             +(1+\kappa)^2\{ay-(1-a)\kappa\}^2\bigr].\nonumber
\end{align}
The leading loss
\(L_0(a)=(1-a)^2\kappa^2+a^2y\) separates the shrinkage bias from
sampling variability. Its minimizer is
\[
 a_* = \frac{\kappa^2}{\kappa^2+y},\qquad
 L_0(a_*)=\frac{y\kappa^2}{\kappa^2+y},
\]
obtained by differentiating this quadratic function. The value of \(a_*\)
is determined by the specified population and is treated as fixed.
Formula~\eqref{eq:two-point-variance}
then describes the second-order variation around the corresponding
finite-dimensional loss center. In particular, with \(a=1\), the
leading loss is \(y\) for every \(\kappa\), whereas its limiting
variance is
\[
 4y^2(1+\kappa^2)^2+8y^3(1+\kappa^2),
\]
which increases with spectral dispersion. The leading loss is therefore
the same across these populations, although its fluctuations differ.
\end{example}
\end{revision}

\section{Proofs}\label{sec:proofs}

We prove Theorems~\ref{thm:main} and~\ref{thm:ewss-clt} first, then
the application results in Section~\ref{sec:applications}. The proofs
of Corollary~\ref{cor:sample-loss} and Proposition~\ref{prop:affine-loss}
share a joint Gaussian limit for three polynomial statistics. Their
proofs and the verification of Example~\ref{ex:two-point-shrinkage}
are grouped in Section~\ref{sec:joint-polynomial-proofs}. The supplement
contains the technical estimates used below and the covariance
simplification in Theorem~\ref{thm:weight-t}.

\subsection{Proof of Theorem~\ref{thm:main}}\label{sec:proof-main}

\begin{revision}
We first work with \(z<0\), where the resolvent is uniformly bounded
and its expectation is well defined. For each weight, decompose
\[
 M_n^g(z)=\bigl[\tr\{(\bS_n-zI)^{-1}g(\bSig)\}
             -\E\tr\{(\bS_n-zI)^{-1}g(\bSig)\}\bigr]
             +\E M_n^g(z).
\]
The bias expansion of \citet{HuYangHan2024} gives the limiting mean.
For the random part, a small real perturbation preserves positive
definiteness, allowing us to apply the joint LSS theorem of
\citet{LiLiYao2018} to the perturbed log determinants.
A uniform finite-difference estimate transfers that joint limit to the
weighted resolvents. Finally, local tightness and uniqueness of analytic
continuation give convergence on the integration contours.

We use two estimates, stated below as Lemmas~\ref{lem:marginal-input}
and~\ref{lem:logdet-bound}. Their proofs and the auxiliary
log-determinant calculations are given in Section~\ref{supp-supp-sec:main} of the supplement.
\end{revision}

We begin with the mean and tightness estimates. For an admissible weight \(g\), write \(u_n=\um_{0,n}(z)\) and set
\begin{align*}
 P_{g,n}(z)&=y_n\int\frac{tg(t)}{(1+tu_n)^2}\,dH_p(t),
 &Q_{g,n}(z)&=y_n\int\frac{t^2g(t)}{(1+tu_n)^3}\,dH_p(t),\\
 R_n(z)&=y_n\int\frac{t^2}{(1+tu_n)^3}\,dH_p(t),
 &\Delta_n(z)&=1-y_n\int\frac{t^2u_n^2}{(1+tu_n)^2}\,dH_p(t).
\end{align*}
We use the localized process and the space \(\mathcal H(\Omega)\)
defined in Section~\ref{sec:preliminaries}.

\begin{lemma}\label{lem:marginal-input}
Under \((H_1)\)--\((H_4)\), the processes
\(M_{n,\mathrm{loc}}^g\) are tight in \(\mathcal H(\Omega)\).
For every compact \(J\subset\Omega\),
\[
 \sup_n\E\sup_{z\in J}|M_{n,\mathrm{loc}}^g(z)|^2<\infty.
\]
For \(z\) in a compact interval of \((-\infty,0)\),
\begin{equation*}
 \E M_n^g(z)=\frac{u_n^2}{z\Delta_n(z)}
       \left\{\frac{P_{g,n}(z)u_nR_n(z)}{\Delta_n(z)}
                      -Q_{g,n}(z)\right\}+o(1),
\end{equation*}
uniformly on that interval. The deterministic expression converges to
\(\mu_g(z)\) in \eqref{eq:lemMean}. At each negative argument,
\(\sup_n\operatorname{Var}\{\tr((\bS_n-zI)^{-1}g(\bSig))\}<\infty\).
\end{lemma}

Lemma~\ref{lem:marginal-input} supplies the limiting mean and the
local moment bounds. Its proof is given in Section~\ref{supp-supp-sec:marginal}. The bias
uses the expansion of \citet{HuYangHan2024}, in the resolvent form
of Lemma C.1 and Section D.1.5 of \citet{HuYangHanPreprint2024}.
The tightness estimate is proved directly and applies throughout
\(\Omega\), including the gap between zero and the nonzero spectrum
when \(y>1\).

We next prove joint convergence for negative arguments.
\rev{The identity
\[
 \left.\partial_s\log\det(\bS_n-zI+s\,g(\bSig))\right|_{s=0}
       =\tr\{(\bS_n-zI)^{-1}g(\bSig)\}
\]
suggests deriving the joint weighted-resolvent limit from a joint
log-determinant limit. We must also differentiate the matching
deterministic center and control the finite-difference error uniformly
in \(n\).} Fix \(z<0\), and introduce a scalar perturbation \(s\):
\begin{align}
 G&=g(\bSig),\qquad A_g(z,s)=-zI+sG,
 &R_g(z,s)&=(\bS_n+A_g(z,s))^{-1},\notag\\
 b_{g,n}(z,s)^{-1}
 &=1+\frac1n\tr\{\bSig(A_g(z,s)+b_{g,n}(z,s)\bSig)^{-1}\}.
 \label{eq:deformed-finite}
\end{align}
For small real \(s\), \(A_g(z,s)\) is positive definite.
The positive solution of \eqref{eq:deformed-finite} is unique and
bounded away from zero. Indeed, multiplying that equation by
\(b=b_{g,n}\) gives
\[
 b+\frac1n\sum_{j=1}^p\frac{bt_j}{-z+sg(t_j)+bt_j}=1,
\]
whose left side is strictly increasing in \(b>0\).
The stability estimate in Section~\ref{supp-supp-sec:logdet} shows that the solution admits
an analytic continuation to a fixed small complex neighborhood of
zero. This neighborhood can be chosen simultaneously for any fixed
finite collection of pairs \((z,g)\).

Put \(T_{g,n}=(A_g+b_{g,n}\bSig)^{-1}\), and define
\begin{equation*}
 \mathcal L_{0,n}^g(z,s)
 =\log\det(A_g+b_{g,n}\bSig)-n\log b_{g,n}-n(1-b_{g,n}).
\end{equation*}
The logarithms are analytic branches continued from positive real
matrices. The coefficient of \(\partial_sb_{g,n}\) in its derivative
is \(\tr(\bSig T_{g,n})-n/b_{g,n}+n=0\). Hence
\[
 \partial_s\mathcal L_{0,n}^g(z,s)=\tr\{T_{g,n}(z,s)G\}.
\]
At zero, \(b_{g,n}(z,0)=-z\um_{0,n}(z)\), so
\begin{equation}\label{eq:logdet-derivative}
 \left.\partial_s\{
 \log\det(\bS_n+A_g(z,s))-\mathcal L_{0,n}^g(z,s)\}
 \right|_{s=0}=M_n^g(z).
\end{equation}

For real \(s\), set
\[
 \mathcal Q_{g,n}=A_g^{-1/2}\bSig^{1/2},\qquad
 B_{g,n}=n^{-1}\mathcal Q_{g,n}\bZ\bZ^*\mathcal Q_{g,n}^*.
\]
Then
\(\log\det(\bS_n+A_g)=\log\det A_g+\tr\log(I+B_{g,n})\).
\rev{The transformed population covariance is
\(\mathcal T_{g,n}=\mathcal Q_{g,n}\mathcal Q_{g,n}^*
=A_g^{-1/2}\bSig A_g^{-1/2}\), with eigenvalues
\(t_j/\{-z+sg(t_j)\}\).} For any fixed collection of perturbations,
the populations are uniformly positive and bounded, their covariance
matrices commute, and their joint empirical spectral measures
converge. The last assertion follows because all the transformed
populations are functions of the same \(t_j\), with discontinuities
on an \(H\)-null set. The entry conditions in
\citet[Theorem 2.1]{LiLiYao2018} hold with \(\alpha_x=1\) and
\(\beta_x=0\). Thus its joint LSS theorem applies to
\(f(x)=\log(1+x)\).

The deterministic LSS center is
\(\mathcal L_{0,n}^g-\log\det A_g\). To verify this, multiply the
transformed population by \(v\ge0\). Both this expression and the
finite-dimensional MP integral of \(\log(1+x)\) vanish at \(v=0\);
their derivatives are \(n\{1-b(v)\}/v\).
The calculation in Section~\ref{supp-supp-sec:logdet} therefore identifies the centers
exactly.

Let \(b_g=b_g(z,s)\) solve the limiting equation, and define
\begin{equation}\label{eq:deformed-limit}
 b_g^{-1}=1+y\int\frac{t}{-z+sg(t)+b_gt}\,dH(t),
 \qquad
 q_g(z,s;t)=\frac{b_gt}{-z+sg(t)+b_gt}.
\end{equation}
For two pairs of perturbations, write
\[
 a_{g,h}(z_1,s_1;z_2,s_2)
   =y\int q_g(z_1,s_1;t)q_h(z_2,s_2;t)\,dH(t).
\]
Evaluation of the joint LSS covariance for \(\log(1+x)\) gives
\begin{equation}\label{eq:logdet-covariance}
 \mathcal V_{g,h}(z_1,s_1;z_2,s_2)
       =-2\log\{1-a_{g,h}(z_1,s_1;z_2,s_2)\}.
\end{equation}
For the real perturbations considered here, the logarithm has a
positive argument. The factor two corresponds to real data.
Section~\ref{supp-supp-sec:logdet} evaluates the two spectral integrations and verifies
these facts.

The passage to weighted resolvents uses the following estimate.
\begin{lemma}\label{lem:logdet-bound}
Fix \(z<0\) and a bounded real weight \(g\). Choose \(\rho>0\) so
that \(\operatorname{Re}A_g(z,s)\succeq cI\) for \(|s|\le2\rho\),
with \(c>0\) independent of \(n\). Define
\[
 F_n^g(z,s)=\log\det(\bS_n+A_g(z,s))
                    -\E\log\det(\bS_n+A_g(z,s)).
\]
Then \(F_n^g(z,\cdot)\) is analytic and
\(\sup_n\sup_{|s|\le2\rho}\E|F_n^g(z,s)|^2<\infty\).
For real \(0<|h|<\rho/2\),
\begin{equation}\label{eq:finite-difference-error}
 \sup_n\E\left|
 \partial_sF_n^g(z,0)-\frac{F_n^g(z,h)-F_n^g(z,0)}h
              \right|^2\le Ch^2.
\end{equation}
\end{lemma}

The proof of Lemma~\ref{lem:logdet-bound} is in Section~\ref{supp-supp-sec:logdet-bound}.
We now apply it to the joint log-determinant limit. For fixed real
perturbations, write
\[
 Y_{n,a}=\log\det(\bS_n+A_{g_a}(z_a,s_a))
                    -\mathcal L_{0,n}^{g_a}(z_a,s_a).
\]
The variables \(Y_{n,a}\) use the deterministic centering of the joint
LSS theorem. To center them by their expectations, we first show that
these expectations converge.
Lemma~\ref{lem:logdet-bound} gives
\(\operatorname{Var}(Y_{n,a})\le C\).
Tightness of \(Y_{n,a}\) implies boundedness of \(\E Y_{n,a}\):
otherwise a subsequence with \(|\E Y_{n,a}|\to\infty\), together
with Chebyshev's inequality,
\[
 \mathbb P\{|Y_{n,a}|\le M\}
 \le \frac{C}{(|\E Y_{n,a}|-M)^2},
\]
would contradict tightness. Hence \(\sup_n\E|Y_{n,a}|^2<\infty\),
which implies uniform integrability of \(Y_{n,a}\). Its expectation
therefore converges to the mean of its Gaussian weak limit.
Subtracting these expectations shows that the finite collection of
\(F_n^{g_a}(z_a,s_a)\) converges to a centered Gaussian vector with
covariance \eqref{eq:logdet-covariance}.

For a fixed nonzero real \(h\), define
\[
 D_{n,a}(h)=\frac{F_n^{g_a}(z_a,h)-F_n^{g_a}(z_a,0)}h,
 \qquad U_{n,a}=\partial_sF_n^{g_a}(z_a,0).
\]
Apply the preceding joint theorem to all \(2r\) triples used in
these quotients. Then \(D_n(h)\) has a centered Gaussian limit,
say \(D(h)\), with covariance
\begin{align*}
 \operatorname{Cov}\{D_a(h),D_b(h)\}
 =h^{-2}\{&\mathcal V_{g_a,g_b}(z_a,h;z_b,h)
       -\mathcal V_{g_a,g_b}(z_a,h;z_b,0)\\
       &-\mathcal V_{g_a,g_b}(z_a,0;z_b,h)
       +\mathcal V_{g_a,g_b}(z_a,0;z_b,0)\}.
\end{align*}
Analyticity of \(\mathcal V\) implies that these covariance matrices
converge as \(h\to0\), with entries
\[
 \left.\partial_{s_1}\partial_{s_2}
       \mathcal V_{g_a,g_b}(z_a,s_1;z_b,s_2)\right|_{s_1=s_2=0}.
\]
Thus \(D(h)\) converges to a centered Gaussian vector \(U\), possibly
with a singular covariance matrix.

For \(t\in\mathbb R^r\), \eqref{eq:finite-difference-error} gives
\begin{align*}
 \left|\E e^{it^{\mathsf T}U_n}
               -\E e^{it^{\mathsf T}D_n(h)}\right|
 &\le\sum_a|t_a|\,\E|U_{n,a}-D_{n,a}(h)|\\
 &\le C_t|h|,
\end{align*}
uniformly in \(n\). First take \(n\to\infty\), and then
\(h\to0\). The characteristic functions converge to that of \(U\),
so the derivative vector converges jointly to \(U\).
By differentiation of the expectation in Lemma~\ref{lem:logdet-bound},
\[
 U_{n,a}=\tr\{(\bS_n-z_aI)^{-1}g_a(\bSig)\}
          -\E\tr\{(\bS_n-z_aI)^{-1}g_a(\bSig)\}.
\]
Thus \(U_n\) is the vector of expectation-centered weighted
resolvents. Lemma~\ref{lem:marginal-input} supplies their deterministic
mean shifts \(\mu_{g_a}(z_a)\).

We now evaluate the covariance of this limiting derivative vector.
Write \(u=\um(z)\), and let \(B_g=\partial_s b_g(z,0)\).
Differentiate \eqref{eq:deformed-limit} at zero, where \(b_g=-zu\):
\[
 B_g\left\{\frac1{z^2u^2}
       -\frac y{z^2}\int\frac{t^2}{(1+tu)^2}\,dH(t)\right\}
 =\frac y{z^2}\int\frac{tg(t)}{(1+tu)^2}\,dH(t).
\]
Since \(u^{-2}-y\int t^2(1+tu)^{-2}dH=1/u'\),
\[
 B_g=yu'\int\frac{tg(t)}{(1+tu)^2}\,dH(t)=d_g(z).
\]
Differentiating the function \(q_g\) gives
\begin{align*}
 \left.\partial_s
       \frac{b_g(z,s)t}{-z+sg(t)+b_g(z,s)t}\right|_{s=0}
 &=\frac{-zt B_g+zu t g(t)}{z^2(1+tu)^2}\\
 &=\frac{t\{u g(t)-d_g(z)\}}{z(1+tu)^2}=w_g(z,t).
\end{align*}
For
\(a=a_{g,h}(z_1,s_1;z_2,s_2)\), it follows that
\begin{align*}
 a|_0&=y\int q_1q_2\,dH,\\
 \partial_{s_1}a|_0&=y\int w_g(z_1,t)q_2(t)\,dH(t),\\
 \partial_{s_2}a|_0&=y\int q_1(t)w_h(z_2,t)\,dH(t),\\
 \partial_{s_1}\partial_{s_2}a|_0
       &=y\int w_g(z_1,t)w_h(z_2,t)\,dH(t).
\end{align*}
Finally,
\[
 \left.\partial_{s_1}\partial_{s_2}[-2\log(1-a)]\right|_0
 =2\left.\left\{\frac{\partial_{s_1}\partial_{s_2}a}{1-a}
       +\frac{(\partial_{s_1}a)(\partial_{s_2}a)}{(1-a)^2}
          \right\}\right|_0.
\]
Substitution gives the cross-covariance function
\eqref{eq:cross-kernel}.

It remains to extend the convergence from the negative real axis to
\(\Omega\). Lemma~\ref{lem:marginal-input} gives tightness of every
localized coordinate \(M_{n,\mathrm{loc}}^{g_a}\) in
\(\mathcal H(\Omega)\). For a fixed finite
collection, choose compact sets for the individual coordinates with
complement probabilities summing to less than \(\varepsilon\).
Their finite product is compact, so the vector of analytic processes
is tight in \(\mathcal H(\Omega)^r\).

Take any convergent subsequence, with limit \(\mathcal M\). Let
\(D\) be the rational points of a fixed negative interval. The
finite-dimensional joint limits on \(D\) have been determined in the
preceding argument, since localization changes these variables only on an event whose
probability tends to zero.
The domain \(\Omega\) is connected. By the identity theorem,
evaluation on \(D\) is injective on \(\mathcal H(\Omega)\).
This function space is separable and complete, and its Borel law is
determined by the evaluations on \(D\). Hence the law of
\(\mathcal M\) is the same for every convergent subsequence.

We verify that evaluations at arbitrary complex arguments are jointly
Gaussian. Start at a rational point in the negative interval. Its
Taylor coefficients are almost sure limits of finite differences of
evaluations on \(D\). Each finite collection of these differences
is jointly Gaussian. Taylor polynomials converge to the analytic
function on a disk centered at that point; their coefficients and
values are therefore jointly Gaussian as well. A finite chain of
overlapping disks reaches any prescribed finite collection of points
in \(\Omega\). Repeating the argument along this chain expresses
the desired evaluations as successive limits of linear combinations
of the initial Gaussian evaluations. The class of finite-dimensional
Gaussian distributions is closed under weak convergence. Thus
\(\mathcal M\) is a Gaussian analytic process.

The local supremum bounds in Lemma~\ref{lem:marginal-input} pass
to the limit by truncation of the squared supremum. They imply local integrability of
\(\sup|\mathcal M|^2\). Cauchy's formula and Fubini's theorem
then show that the mean is analytic and the bilinear covariance is
analytic in its two arguments. Their values on the negative axis
are \(\mu_g\) and \(\Gamma_{g,h}\); successive application of the
identity theorem gives these formulas throughout \(\Omega\).
For real weights,
\(M_{n,\mathrm{loc}}^g(\overline z)
      =\overline{M_{n,\mathrm{loc}}^g(z)}\), and the limit has the
same symmetry. Its Hermitian covariance is consequently
\(\Gamma_{g,h}(z_1,\overline{z_2})\).
Every subsequential limit is now identified. This proves convergence
of the full sequence in \(\mathcal H(\Omega)^r\). Restriction to
the chosen finite union of contours is continuous, completing the
proof of Theorem~\ref{thm:main}.

\subsection{Proof of Theorem~\ref{thm:ewss-clt}}

On \(E_n\), subtracting \eqref{eq:weighted-lss-centering} from
\eqref{eq:weighted-lss-cauchy} gives
\[
 p\{T_n(\varphi_a,g_a)-\Theta_n(\varphi_a,g_a)\}
 =-\frac1{2\pi i}\oint_{\mathcal C}
              \varphi_a(z)M_n^{g_a}(z)\,dz.
\]
For a fixed piecewise smooth contour, the integration map is continuous
under the uniform norm: its absolute value is at most
\((2\pi)^{-1}\operatorname{length}(\mathcal C)
\sup_{\mathcal C}|\varphi_a|\sup_{\mathcal C}|M_n^{g_a}|\).
Theorem~\ref{thm:main} and the continuous mapping theorem therefore
give the joint limit after localization. Since
\(\mathbb P(E_n)\to1\), the unlocalized statistic has the same weak
limit. Integrating the mean and covariance gives
\eqref{eq:lss-mean} and \eqref{eq:lss-cross}. The statistic is real
for real-valued tests on the spectral enclosure, so its limiting vector
is real Gaussian. Zero variances and singular covariance matrices are
allowed.

\subsection{Proof of Theorem~\ref{thm:variance-component}}

Take the contour
\(\mathcal C_+\) surrounding \([L,U]\) and a separate contour
surrounding zero. The function \(\varphi_V\) is \(z^{-1}\) on the
first enclosed neighborhood and zero on the second. On \(E_n\),
functional calculus gives
\[
 V_{\mathsf X,\boldsymbol\beta}
       =\sigma^2y_nT_n(\varphi_V,t).
\]
The tightness estimate in Section~\ref{supp-supp-sec:marginal} applies to
\(\mathcal C_+\), including its real crossing in \((0,L)\). On \(E_n\),
\[
 p\{V_{\mathsf X,\boldsymbol\beta}-\sigma^2y_n\Theta_{V,n}\}
   =\sigma^2y_n\,p\{T_n(\varphi_V,t)-\Theta_n(\varphi_V,t)\}
\]
Theorem~\ref{thm:ewss-clt} and Slutsky's theorem give the stated limit;
\(\mathbb P(E_n)\to1\) removes the localization.

\subsection{Proof of Theorem~\ref{thm:frobenius-loss}}

Let
\(Z_{0,n}=p\{T_n(\ell^2,1)-\Theta_n(\ell^2,1)\}\) and
\(Z_{1,n}=p\{T_n(\ell,t)-\Theta_n(\ell,t)\}\).
Equations~\eqref{eq:frobenius-loss-expansion} and
\eqref{eq:loss-centering} give
\[
 p\{L_n(\ell)-L_{0,n}(\ell)\}=Z_{0,n}-2Z_{1,n}.
\]
Theorem~\ref{thm:ewss-clt} gives a joint limit for
\((Z_{0,n},Z_{1,n})\). The displayed linear map yields
\(\mu_0-2\mu_1\) and \(v_{00}+4v_{11}-4v_{01}\), as required.

\subsection{Polynomial losses and the two-point example}\label{sec:joint-polynomial-proofs}

The sample and affine-shrinkage losses are linear combinations of
\(\tr(\bS_n^2)\), \(\tr(\bS_n\bSig)\), and \(\tr(\bS_n)\).
We first obtain the joint Gaussian limit of these three statistics,
then use the corresponding linear combinations to prove the two
results. Unless explicitly stated otherwise, the covariances below
are parameters of the Gaussian weak limits, rather than assertions
of finite-sample moment convergence.

Define the centered vector by
\begin{align*}
 Z_{A,n}&=\tr(\bS_n^2)-p(H_{2,n}+y_nH_{1,n}^2),\\*
 Z_{B,n}&=\tr(\bS_n\bSig)-pH_{2,n},\qquad
 Z_{C,n}=\tr(\bS_n)-pH_{1,n}.
\end{align*}
Expansion of the finite MP transform at infinity gives
\[
 \Theta_n(x,1)=H_{1,n},\qquad
 \Theta_n(x,t)=H_{2,n},\qquad
 \Theta_n(x^2,1)=H_{2,n}+y_nH_{1,n}^2.
\]
For example, substituting
\(\um_{0,n}(z)=-z^{-1}-y_nH_{1,n}z^{-2}+O(z^{-3})\)
into the finite-dimensional weighted transform in \eqref{eq:finite-centering}
gives the first two weighted moments. The limiting mean has expansion
\[
 \mu_g(z)=-y\int t^2g(t)\,dH(t)\,z^{-3}+O(z^{-4}).
\]
The coefficient of \(z^{-3}\), with the minus sign in the Cauchy formula,
therefore gives \(\E Z_A=yH_2\); the other two limiting means are zero.

Theorem~\ref{thm:ewss-clt} now gives a joint Gaussian limit
\((Z_A,Z_B,Z_C)\).

We calculate the covariance entries, starting with the linear
spectral tests. These calculations also give
Corollary~\ref{col:linear-tests}.

For \(\varphi=1\),
\[
 T_n(1,g)=p^{-1}\tr g(\bSig)=\Theta_n(1,g)
\]
exactly. The scalar limiting mean and variance are therefore zero.
Equivalently, \(w_g(z,t)=O(z^{-2})\) and
\(\mu_g(z)=O(z^{-3})\) at infinity, so the contour integrals for a
constant test vanish.
For \(\varphi(x)=x\), independence of the normalized columns and
\eqref{supp-eq:quadratic-identity} give
\begin{align}
 \E\tr\{\bS_ng(\bSig)\}&=\tr\{\bSig g(\bSig)\},\notag\\*
 \operatorname{Cov}\bigl(\tr\{\bS_ng(\bSig)\},
                         \tr\{\bS_nh(\bSig)\}\bigr)
 &=\frac2n\tr\{\bSig^2g(\bSig)h(\bSig)\}.
 \label{eq:linear-test-cov}
\end{align}
The limiting bias is zero, and the limiting covariance is
\(2y\int t^2g(t)h(t)\,dH(t)\). This is also the
\(z_1^{-2}z_2^{-2}\) coefficient of \eqref{eq:cross-kernel}.

For the covariance between a quadratic and a linear statistic,
we expand the general resolvent covariance at infinity.

Write \(\langle f\rangle=\int f(t)\,dH(t)\),
\(H_k=\langle t^k\rangle\), and \(v=z^{-1}\). Expansion at
infinity gives
\begin{align*}
 \um(z)&=-v-yH_1v^2+O(v^3),\\*
 q(z,t)&=-tv-(t^2+yH_1t)v^2+O(v^3),\\*
 w_g(z,t)&=-tg(t)v^2
       -t\{2tg(t)+yH_1g(t)+y\langle tg\rangle\}v^3+O(v^4).
\end{align*}
The second term in \eqref{eq:cross-kernel} starts at order
\(z_1^{-3}z_2^{-3}\). Thus the \(z_1^{-3}z_2^{-2}\)
coefficient is obtained by multiplying the third-order term of
\(w_g(z_1,t)\) and the second-order term of \(w_h(z_2,t)\):
\begin{equation}\label{eq:mixed-polynomial-coefficient}
 4y\langle t^3gh\rangle
    +2y^2\{H_1\langle t^2gh\rangle
                      +\langle tg\rangle\langle t^2h\rangle\}.
\end{equation}
This coefficient depends on the pairing of the weights with the two
spectral tests.

The linear covariance formula
\eqref{eq:linear-test-cov} yields
\[
 V_{BB}=2yH_4,\qquad V_{BC}=2yH_3,\qquad V_{CC}=2yH_2.
\]
For the mixed terms, set \(g=1,h=t\) and \(g=h=1\), respectively,
in \eqref{eq:mixed-polynomial-coefficient}. This gives
\[
 V_{AB}=4yH_4+4y^2H_1H_3,\qquad
 V_{AC}=4yH_3+4y^2H_1H_2.
\]
It remains to extract the \(z_1^{-3}z_2^{-3}\) coefficient for
\(V_{AA}\). For \(g=1\),
\[
 w_1(z,t)=-tz^{-2}-2t(t+yH_1)z^{-3}+O(z^{-4}).
\]
The product of the two third-order terms in the first term of
\eqref{eq:cross-kernel} contributes
\(8yH_4+16y^2H_1H_3+8y^3H_1^2H_2\).
Also,
\(D_{12}^{-1}=1+yH_2z_1^{-1}z_2^{-1}+
O(z_1^{-2}z_2^{-1})+O(z_1^{-1}z_2^{-2})\)
as a bivariate expansion at infinity. Its first correction contributes
\(2y^2H_2^2\). The leading term of the product of the two integrals
in the second term of \eqref{eq:cross-kernel} contributes another
\(2y^2H_2^2\). Thus
\begin{equation*}
 V_{AA}=8yH_4+16y^2H_1H_3+4y^2H_2^2+8y^3H_1^2H_2.
\end{equation*}
Writing \(V_{AB}=\operatorname{Cov}(Z_A,Z_B)\), and similarly
for the other entries, we have obtained
\begin{align*}
 V_{AA}&=8yH_4+16y^2H_1H_3+4y^2H_2^2+8y^3H_1^2H_2,\\*
 V_{AB}&=4yH_4+4y^2H_1H_3,&
 V_{AC}&=4yH_3+4y^2H_1H_2,\\*
 V_{BB}&=2yH_4,& V_{BC}&=2yH_3,\qquad V_{CC}=2yH_2.
\end{align*}

The joint limit is now fully specified and can be applied to both losses.

\paragraph{Proof of Corollary~\ref{cor:sample-loss}.}

For the sample covariance, expanding \(\|\bS_n-\bSig\|_F^2\) and
using the finite MP moments gives
\[
 p\{p^{-1}\|\bS_n-\bSig\|_F^2-y_nH_{1,n}^2\}
     =Z_{A,n}-2Z_{B,n}.
\]
Its limiting mean is \(yH_2\), and its variance is
\[
 V_{AA}+4V_{BB}-4V_{AB}=4y^2H_2^2+8y^3H_1^2H_2.
\]
This proves Corollary~\ref{cor:sample-loss}.
For the finite-sample expectation in the main text, expand
\(\tr(\bS_n^2)=n^{-2}\sum_{i,j}(\bx_i^*\bx_j)^2\).
When \(i\ne j\), its summand has expectation \(\tr\bSig^2\).
When \(i=j\), the fourth-moment identity gives
\(\E(\bx_i^*\bx_i)^2=(\tr\bSig)^2+2\tr\bSig^2\).
Consequently,
\[
 \E\tr(\bS_n^2)=\left(1+\frac1n\right)\tr\bSig^2
                          +\frac1n(\tr\bSig)^2.
\]
Subtract \(2\E\tr(\bS_n\bSig)-\tr\bSig^2=\tr\bSig^2\)
and divide by \(p\) to obtain the stated expectation.

\paragraph{Proof of Proposition~\ref{prop:affine-loss}.}

For the general fixed affine rule,
direct expansion gives
\begin{align*}
 p^{-1}\|a\bS_n+bI-\bSig\|_F^2
 &=a^2p^{-1}\tr(\bS_n^2)-2ap^{-1}\tr(\bS_n\bSig)
            +2abp^{-1}\tr(\bS_n)\\*
 &\quad+H_{2,n}-2bH_{1,n}+b^2.
\end{align*}
Insert the three finite-dimensional centers. Their deterministic terms
are exactly \(L_{0,n}(a,b)\) in \eqref{eq:affine-center}; subtracting
them gives
\begin{equation*}
 p\{L_n(\ell_{a,b})-L_{0,n}(a,b)\}
          =a^2Z_{A,n}-2aZ_{B,n}+2abZ_{C,n}.
\end{equation*}
Hence the limit is \(a^2Z_A-2aZ_B+2abZ_C\), with mean
\(a^2yH_2\).

Set \(U=Z_A-2Z_B-2yH_1Z_C\). The covariance calculations above imply
\begin{align*}
 \operatorname{Cov}(U,Z_B)
   &=(4yH_4+4y^2H_1H_3)-4yH_4-4y^2H_1H_3=0,\\*
 \operatorname{Cov}(U,Z_C)
   &=(4yH_3+4y^2H_1H_2)-4yH_3-4y^2H_1H_2=0.
\end{align*}
Substituting all six covariance entries into
\[
 \operatorname{Var}(U)=V_{AA}+4V_{BB}+4y^2H_1^2V_{CC}
                 -4V_{AB}-4yH_1V_{AC}+8yH_1V_{BC}
\]
leaves \(4y^2H_2^2\). In particular, \(U\) is uncorrelated with every
linear combination of \(Z_B,Z_C\). The covariance of two such
combinations is
\[
 \operatorname{Cov}(\alpha Z_B+\beta Z_C,\gamma Z_B+\delta Z_C)
       =2y\int t^2(\alpha t+\beta)(\gamma t+\delta)\,dH(t).
\]
Write the limiting affine loss as
\[
 a^2U+2a\{(a-1)Z_B+(b+ayH_1)Z_C\}.
\]
Its variance is \eqref{eq:affine-variance}, and its covariance with
the corresponding expression for \((c,d)\) is \eqref{eq:affine-cross}.
Every finite collection is a linear image of the same Gaussian vector,
which proves the joint assertion of Proposition~\ref{prop:affine-loss}.

\paragraph{Proof of Example~\ref{ex:two-point-shrinkage}.}
Substitute the two-point population into the formulas of
Proposition~\ref{prop:affine-loss}. The finite population and its
limit both assign mass \(1/2\) to \(1-\kappa\) and \(1+\kappa\).
Set \(b=1-a\) in \eqref{eq:affine-center}, use
\(H_{1,n}=1,H_{2,n}=1+\kappa^2\), and collect terms to obtain
\((1-a)^2\kappa^2+a^2y_n\).
For the variance integrand, the two values of
\((a-1)t+1-a+ay\) are
\(ay+(1-a)\kappa\) and \(ay-(1-a)\kappa\), respectively.
Averaging their squared values with weights \(t^2/2\) gives
\eqref{eq:two-point-variance}. Finally,
\[
 L_0'(a)=2\{(\kappa^2+y)a-\kappa^2\},\qquad
 L_0''(a)=2(\kappa^2+y)>0.
\]
Thus \(a_*\in(0,1)\) is the unique minimizer on \([0,1]\), with
value stated in the example. The CLT remains centered at the
finite-dimensional expression containing \(y_n\).

\begingroup
\setstretch{1.1}
\interlinepenalty=10000
\small
\bibliographystyle{plainnat}
\bibliography{ref}
\endgroup
\end{document}